\documentclass{amsart} 
\usepackage{amssymb,latexsym,amsfonts,verbatim,amscd,ifthen,blkarray} 
\usepackage{color}
\usepackage[enableskew]{youngtab}

\numberwithin{equation}{section}

\theoremstyle{plain}

\newtheorem{theorem}[equation]{Theorem}   
\newtheorem{lemma}[equation]{Lemma} 
\newtheorem{proposition}[equation]{Proposition} 
\newtheorem{corollary}[equation]{Corollary} 
 
\newtheorem{thm}[equation]{Theorem}
\newtheorem{prop}[equation]{Proposition}

\theoremstyle{definition}

\newtheorem{definition}[equation]{Definition} 
\newtheorem{remark}[equation]{Remark} 
\newtheorem{remarks}[equation]{Remarks}
 
\newtheorem{Example}[equation]{Example} 
\newtheorem{example}[equation]{Example} 
\newtheorem{dfn}[equation]{Definition}

\DeclareMathOperator{\HH}{H}

\DeclareMathOperator{\id}{id}  
\DeclareMathOperator{\coker}{Coker} 
\DeclareMathOperator{\Coker}{Coker} 
 
\DeclareMathOperator{\grade}{grade} 
\DeclareMathOperator{\rank}{rank}
\DeclareMathOperator{\im}{Im} 
\DeclareMathOperator{\pd}{pd}

\DeclareMathOperator{\Spec}{Spec}

\DeclareMathOperator{\start}{start}
\DeclareMathOperator{\finish}{end}

\begin{document}   

\renewcommand{\:}{\! :} 
\newcommand{\p}{\mathfrak p} 
\newcommand{\m}{\mathfrak m}
\newcommand{\e}{\epsilon}
\newcommand{\lra}{\longrightarrow} 
\newcommand{\ra}{\rightarrow} 
\newcommand{\altref}[1]{{\upshape(\ref{#1})}} 
\newcommand{\bfa}{\boldsymbol{\alpha}} 
\newcommand{\bfb}{\boldsymbol{\beta}} 
\newcommand{\bfg}{\boldsymbol{\gamma}} 
\newcommand{\bfM}{\mathbf M} 
\newcommand{\bfI}{\mathbf I} 
\newcommand{\bfC}{\mathbf C} 
\newcommand{\bfB}{\mathbf B} 
\newcommand{\bsfC}{\bold{\mathsf C}} 
\newcommand{\bsfT}{\bold{\mathsf T}}
\newcommand{\smsm}{\smallsetminus} 
\newcommand{\ol}{\overline}

\newlength{\wdtha}
\newlength{\wdthb}
\newlength{\wdthc}
\newlength{\wdthd}
\newcommand{\elabel}[1]
           {\label{#1}  
            \setlength{\wdtha}{.4\marginparwidth}
            \settowidth{\wdthb}{\tt\small{#1}} 
            \addtolength{\wdthb}{\wdtha}
            \raisebox{\baselineskip}
            {\color{red} 
             \hspace*{-\wdthb}\tt\small{#1}\hspace{\wdtha}}}  

\newcommand{\mlabel}[1] 
           {\label{#1} 
            \setlength{\wdtha}{\textwidth}
            \setlength{\wdthb}{\wdtha} 
            \addtolength{\wdthb}{\marginparsep} 
            \addtolength{\wdthb}{\marginparwidth}
            \setlength{\wdthc}{\marginparwidth}
            \setlength{\wdthd}{\marginparsep}
            \addtolength{\wdtha}{2\wdthc}
            \addtolength{\wdtha}{2\marginparsep} 
            \setlength{\marginparwidth}{\wdtha}
            \setlength{\marginparsep}{-\wdthb} 
            \setlength{\wdtha}{\wdthc} 
            \addtolength{\wdtha}{1.4ex} 
            \settowidth{\wdthb}{\tt\small{#1}} 
            \marginpar{\vspace*{\baselineskip}
\smash{\raisebox{0.7\baselineskip}{\tt\small{#1}}\hspace{-\wdthb}%
\raisebox{.3\baselineskip}
            {\rule{\wdtha}{0.5pt}}} }
            \setlength{\marginparwidth}{\wdthc} 
            \setlength{\marginparsep}{\wdthd}  }  
            


\title{Torsion freeness of Schur modules}
\author[M. Allahverdi]{Muberra Allahverdi} 
\address{Department of Mathematics\\
         University at Albany, SUNY\\ 
         Albany, NY 12222}
\email{mallahverdi@albany.edu}

\author[A. Tchernev]{Alexandre Tchernev} 
\address{Department of Mathematics\\
         University at Albany, SUNY\\ 
         Albany, NY 12222}
\email{atchernev@albany.edu}


\begin{abstract}
Let $R$ be a Noetherian commutative ring and $M$ a 
$R$-module with $\pd_R M\le 1$ that has rank.
Necessary and sufficient conditions were provided 
in \cite{L1974} for an exterior power 
$\wedge^k M$ to be torsion free. When $M$ is an 
ideal of $R$ similar necessary and sufficient 
conditions were provided 
in \cite{T2007} for a symmetric power $S_k M$  
to be torsion free. We extend these results to a 
broad class of Schur modules $L_{\lambda/\mu} M$. 
En route, for any map 
of finite free $R$ modules 
$\phi\: F\rightarrow G$ 
we also study the general structure of the  
Schur complexes $L_{\lambda/\mu}\phi$, 
and provide 
necessary and sufficient conditions for the 
acyclicity of any given $L_{\lambda/\mu}\phi$
by computing explicitly the 
radicals of the ideals of maximal minors of 
all its differentials. 
\end{abstract}

\maketitle

\section*{Introduction}

Torsion freeness of symmetric powers of ideals, and more generally of modules, has been studied extensively in commutative algebra and algebraic geometry, see for example \cite{VW1994} and the references there. 
When $R$ is a commutative Noetherian ring and $k \geq 2$, equivalent conditions were provided in \cite{T2007} for the $k$th symmetric power $S_k I$ of an ideal $I$ to be torsion free when the projective dimension of $I$ is less than or equal to 1. 
In a similar vein, due to their significance for 
the study of the structure of finite free 
resolutions, see \cite{BE1974}, torsion freeness 
of exterior powers of $R$-modules was 
investigated in \cite{L1974, L1975, TW2004}, and 
similar necessary and sufficient conditions for 
torsion freeness were provided in \cite{L1974} 
for a module $M$ of projective dimension 
at most $1$ that has rank. 

Under the same assumptions on the 
module $M$, 
in this paper we investigate more 
generally the torsion freeness of its Schur 
modules $L_{\lambda/\mu}M$, 
see 
Section~\ref{S:Schur-Weyl} for the 
definitions. In our second main result, 
Theorem~\ref{T:second-main}, we provide, 
for a broad family of skew shapes 
$\lambda/\mu$,  a necessary 
and sufficient condition for the torsion 
freeness of a given Schur module 
$L_{\lambda/\mu}M$ in 
terms of the grades of the Fitting ideals of 
$M$. This generalizes results from 
\cite{L1974,T2007} and also from the first
author's thesis \cite{A2017}. 

To obtain this result it was necessary, 
for any given map $\phi\: F\longrightarrow G$ 
of finite free $R$-modules, to 
investigate in detail the structure of the 
Schur complexes $L_{\lambda/\mu}\phi$ 
introduced in \cite{ABW1982}. In our 
first main result, Theorem~\ref{sqrt}, and 
its corollary Theorem~\ref{acyclicity}, 
we compute explicitly the radicals of the 
ideals of maximal minors of the differentials 
of $L_{\lambda/\mu}\phi$ and provide a 
necessary and sufficient condition for its 
acyclicity. These results are of independent 
interest, and generalize work from    
\cite{ABW1982} and \cite{MS1996, MS1999}. 
We also describe completely 
in Theorem~\ref{genericacyclicity} the 
skew shapes $\lambda/\mu$ for which   
$L_{\lambda/\mu}\phi$ is acyclic in the 
generic case. 

The paper is organized as follows. In 
Section~\ref{intro.ch} we recall basic facts 
from commutative algebra, and establish 
notation. In Section~\ref{S:rigidity} 
we review the 
theory of rigid functors from \cite{T2007}.  
This will be one of the main tools we use to 
study torsion freeness. 
Section~\ref{S:Schur-Weyl} is 
devoted to recalling 
definitions and basic properties of 
Schur and Weyl modules that we will need. 
In Section~\ref{prelimsymmetric} 
we do the same for 
Schur complexes. In 
Section~\ref{S:threshhold} we 
introduce the notion of \emph{threshold 
number} and describe some of its 
elementary properties. In 
Section~\ref{prob.ch} 
we study in detail the structure of the 
differentials of the Schur complexes. 
We state and prove our first main result, 
Theorem~\ref{sqrt}, and obtain important 
corollaries. Finally, in 
Section~\ref{S:torsion-freeness} 
we use these results to provide in  
Theorem~\ref{T:second-main} 
a necessary and sufficient condition for 
the torsion freeness of a Schur module.

\section{Preliminaries}
\label{intro.ch}

Throughout this paper 
rings are commutative Noetherian with unit,
modules are unitary and finitely generated, 
and chain complexes are zero in negative 
homological degrees.

Let $\mathbb F=(F_i,\phi_i)$ be a finite (chain) complex 
of free $R$-modules of finite rank.
If $\mathbb F$ is nonzero, we set 
\begin{align*}
\start \mathbb F &= \min\{n\mid F_n\ne 0\} \\
\finish\mathbb F &= \max\{n\mid F_n\ne 0\}. 
\end{align*}
The \emph{expected rank} of the differential 
$\phi_n$ is the integer 
\[
r_n= \sum_{i\ge n} (-1)^{i-n}\rank F_i. 
\]
We write $I(\phi_n)$ for the ideal 
$I_{r_n}(\phi_n)$ generated by all minors of 
size $r_n$ of $\phi_n$. We say that $\mathbb F$ 
is \emph{acyclic} if $\HH_n\mathbb F=0$ for 
$n\ne 0$, and call $\mathbb F$ \emph{exact} if 
in addition $\HH_0\mathbb F=0$.  

When 
$
F\xrightarrow{\phi} G\rightarrow M\rightarrow 0
$
is a finite free presentation of a 
finitely generated $R$-module $M$ 
with $\rank G=g$, the $k$th \emph{Fitting ideal} 
$\operatorname{Fitt}_k(M)$ of $M$ is the ideal 
$I_{g-k}(\phi)$. It is well known, see e.g. 
\cite[Section 20.2]{E1995}, that this is an invariant of $M$ 
and does not depend on the choice of the finite 
free presentation. 

We recall some of the acyclicity results of 
Buchsbaum and Eisenbud \cite{BE1973}.  

\begin{thm}[Buchsbaum-Eisenbud Acyclicity 
Criterion] \label{eisenbudacyclicity}
Let $R$ be a ring. A complex
\[
\mathbb{F}: 
0 \rightarrow F_n \xrightarrow{\phi_n} F_{n-1}   
\rightarrow \cdots \rightarrow F_1 
\xrightarrow{\phi_1} F_0 \rightarrow 0 
\]
of finite free $R$-modules is acyclic  
if and only if 

\[
\operatorname{grade}I(\phi_k) \geq k
\]
for all $k\ge 1$. In that case, we also have 
\[
\sqrt{I(\phi_1)}\subseteq \dots \subseteq 
\sqrt{I(\phi_i)}\subseteq \dots 
\]
\end{thm}

We also need the following basic consequence of
the Buchsbaum-Eisenbud acyclicity criterion. 

\begin{prop}\label{torsionfreeness}
Let $\mathbb{F} = (F_i,\phi_i)$ be a free 
resolution of a module $M$. The following are 
equivalent:
\begin{enumerate}
\item
M is torsion free;
\item
$\mathbb{F}/r\mathbb{F}$ is a resolution of 
$M/rM$ over $R/(r)$ for every nonzero divisor 
$r \in R$;
\end{enumerate}
If $\mathbb{F}$ is finite, conditions $(1)$ and 
$(2)$ are also equivalent to 
\begin{enumerate}
\item[(3)]
$\operatorname{grade} I(\phi_i) \geq i+1$ 
for each $i \geq 1$.  
\end{enumerate}
\end{prop}

In the coming sections, 
for a map of finite free modules 
$\phi\: F \lra G$ with $f=\rank F<\rank G$, 
we will need to use 
a certain map $\epsilon$ coming from the 
so-called Buchsbaum-Rim complex of $\phi$,  
see \cite{E1995}. 
The map 
\begin{equation}\label{E:epsilon-def}
\epsilon=\epsilon(\phi): 
\wedge^{f+1}G^* \otimes \wedge^f F 
\rightarrow G^*
\end{equation}  
is defined as the following composition: 
{\small
\[
\wedge^{f+1}G^*\otimes\wedge^f F 
\xrightarrow{\Delta\otimes 1}
G^*\otimes\wedge^f G^*\otimes\wedge^f F 
\xrightarrow{1\otimes\wedge\phi^*\otimes 1}
G^*\otimes\wedge^f F^*\otimes\wedge^f F 
\xrightarrow{1\otimes\mu} 
G^*\otimes R = G^*, 
\]
}
\negthickspace
where $\Delta$ is the diagonal map, and 
$\mu$ is the evaluation map. 
In terms of elements, we have the following 
formula: 
\begin{equation}\label{E:epsilon}
\epsilon(e^*_I \otimes b) = 
\sum_{J \subset I, |J|=f} 
\operatorname{sgn}(J \subset I)
(\operatorname{det} \phi_J)e^*_{I \backslash J}
\end{equation}
where $I = \{i_1, \dots, i_{f+1}\}$ is a subset 
of $\{1, \dots, g\}$, 
the set $\{e_1,\dots, e_g\}$ is a basis of $G$, 
the set $\{e_1^*,\dots, e_g^*\}$ is the 
corresponding dual basis of $G^*$, 
the element 
$e_I = e_{i_1} \wedge \cdots \wedge e_{i_{f+1}}$ 
is a basis element of $\wedge^{f+1} G^*$, 
the set $\{b_1,\dots, b_f\}$ is a basis of $F$, 
$b=b_1\wedge\cdots\wedge b_f$ is the free 
generator of $\wedge^f F$, 
the matrix $\phi_J$ is obtained from the 
matrix of $\phi$ 
for the given bases of $F$ and $G$ by taking only 
the columns indexed by elements of $J$, 
and $\operatorname{sgn}(J \subset I)$ is 
the sign of the permutation of $I$ that places 
the elements of $J$ in 
the first $f$ positions.

\begin{remarks}\label{R:epsilon-functorial} 
(a) 
It is routine to check that the composition 
$\phi^*\circ\epsilon=0$, and that, if $F=0$, then 
$\epsilon=\id_{G^*}$. 

(b) 
Given a commutative diagram of free $R$-modules 
\[
\begin{CD} 
F @> \phi >> G   \\ 
@V \psi VV   @VV \gamma V  \\ 
F' @>> \phi' >  G' 
\end{CD}
\]
with $f=\rank F<\rank G=\rank G'$, and 
with $\psi$ an isomorphism, 
it is a routine computation to verify that 
the diagram 
\[
\begin{CD}
\wedge^f F\otimes \wedge^{f+1}G^* 
@>\epsilon(\phi) >> 
G^*                               \\ 
@A \wedge\psi^{-1}\otimes\wedge\gamma^* AA 
@AA \gamma^* A                         \\ 
\wedge^f F'\otimes \wedge^{f+1}{G'}^*
@>>\epsilon(\phi') > 
{G'}^*                               
\end{CD}
\]
is also commutative, where $f=\rank F=\rank F'$. 
\end{remarks}

The following is the key to extending the results 
of \cite{T2007} to the case of Schur modules.

\begin{lemma}\label{exactatg}
Let $R$ be a ring,  
let $M$ be a nonzero torsion-free 
$R$-module with a finite free resolution 
\[
0 \lra F \overset{\phi}\lra G\lra M\lra 0, 
\]
and let $f=\rank F$. Then 
\[
F \overset{\phi}\lra G \xrightarrow{\epsilon^*}
\wedge^f F^*\otimes\wedge^{f+1}G
\]
is exact at $G$, where $\epsilon^*=\epsilon(\phi)^*$ 
is the dual 
of the map $\epsilon(\phi)$ 
from \eqref{E:epsilon-def}.
\end{lemma}

We use the following basic fact, 
see \cite[Theorem 1.3.4]{BH1993}, to prove the 
lemma and reduce to the case where 
$R$ is a field.  

\begin{lemma} \label{lemmab&h}
Let $(R,m,k)$ be a local ring, and 
$\psi: K \rightarrow L$ a homomorphism of 
finitely generated  
$R$-modules. Suppose that $K$ is free, and let 
$M$ be an $R$-module with $m$ an associated 
prime of M. Suppose that $\psi \otimes M$ is 
injective. Then:
\begin{enumerate}
\item 
$\psi \otimes k$ is injective;
    
\item 
if $L$ is a free $R$-module, then $\psi$ is
injective, and $\psi(K)$ is a free direct 
summand of $L$.
\end{enumerate}
\end{lemma}

\begin{proof}[Proof of Lemma \ref{exactatg}]
We want to show that 
$
M = G/(\operatorname{Im}\phi)
\xrightarrow{\overline{\epsilon^*}}  
\wedge^{f+1}G \otimes \wedge^f F^*
$ 
is injective, where 
$\overline{\epsilon^*}$ is the induced natural 
map. 

Since $M$ is torsion-free, it suffices to show 
that 
\smash[t]{
$
M_S \xrightarrow{\overline{\epsilon^*}} 
(\wedge^{f+1} G \otimes \wedge^f F^*)_S
$
}
is injective where 
$
S = 
\{ 
\text{nonzero-divisors \space of \space} R 
\}. 
$ 
Thus we may replace $R$ with $R_S$. Since 
$
\operatorname{Ker}
\bigl(\overline{\epsilon^*}\bigr)_\mathfrak{p} = 0
$ 
$\forall$ $\mathfrak{p}$ implies that 
$
\operatorname{Ker}
\bigl(\overline{\epsilon^*}\bigr) = 0, 
$ 
we may replace $R$ with $R_\mathfrak{p}$. 
Since $M$ has rank, we know that 
$M$ is free over $R$, hence $\phi$ is split. 
Therefore, we may use 
Lemma~\ref{lemmab&h} and it suffices to prove the 
result with $R$ replaced by the field 
$R/\mathfrak{p}R$. 


Since the dimension of 
$\operatorname{Im} \phi = f$, the dimension of 
$\operatorname{Ker} \epsilon^*$ is at least $f$. 
We use the rank-nullity formula 
$
g = \operatorname{dim}G = 
\operatorname{dim}\operatorname{Ker}\epsilon^* + 
\operatorname{dim}\operatorname{Im} \epsilon^*
$ 
to show that it must be exactly $f$.  

 

 Since 
 $
 \operatorname{rank}\epsilon^* = 
 \{ \text{determinantal rank of } \epsilon^*\}, 
 $ 
 it remains to show a $(g-f) \times (g-f)$ 
 submatrix of $\epsilon^*$ with a nonzero 
 determinant. 
 

Since $\phi$ is injective, it must have a nonzero 
$f \times f$ minor. Without loss of generality, 
assume that $d_{1, \dots, f} \not= 0$, where 
$d_{1, \dots, f}$ is the minor corresponding to 
rows $1, \dots, f$. Then, 
with notation as in \eqref{E:epsilon}, the map 
$\epsilon^*$ has the matrix
%
\[
\begin{blockarray}{cccccc}
 \hdots & e_{f+1} & e_{f+2} & \hdots & e_g & \\
 \begin{block}{(ccccc)c}
\hdots & d_{1, \dots, f} & 0 & \hdots & 0 
       & e_{1, \dots, f, f+1} \otimes b^* 
\\
\hdots & * & \ddots & \ddots & \vdots & \vdots
\\
\hdots & \vdots & \ddots & \ddots & 0 
       & e_{1, \dots, f, g-2} \otimes b^*
\\ 
\hdots & * & \hdots & * & d_{1, \dots, f}
       & e_{1, \dots, f, g-1} \otimes b^*
\\
\hdots & * & * & * & * 
       & \ \ \ e_{1, \dots, f, g} \otimes b^*
 \\
\vdots & \vdots & \vdots & \vdots & \vdots 
       & \vdots \\
\end{block} 
\end{blockarray}
\]
and the result follows. 
\end{proof}

\begin{remark}
Suppose that $\phi$ a matrix of indeterminates. 
In this case it follows from \cite{BV} that 
$\operatorname{grade} I_f (\phi) \geq 2$, so that 
$M$ is torsion free by 
Proposition~\ref{torsionfreeness}. Hence, the 
above lemma holds in this generic case.
\end{remark}

\section{Rigidity}
\label{S:rigidity}

We briefly review the notions and results on 
rigid functors from \cite{T2007} 
that will be used in this paper.

%

Let ${\bf Comp}$ be the category with objects $\{R, \mathbb{F} \}$, 
where $R$ is a ring and $\mathbb{F}$ a complex of free $R$-modules, and 
morphisms $\{\rho, \phi\}: \{R, \mathbb{F}\} \rightarrow \{S, \mathbb{G}\}$, 
where $\rho: R \rightarrow S$ is a ring homomorphism and 
$\phi: S \otimes_{\rho} \mathbb{F} \rightarrow \mathbb{G}$ 
is a morphism of complexes over $S$.
Here, $S \otimes_{\rho} \mathbb{F} = S \otimes_R \mathbb{F}$ 
where we consider $S$ as an $R$-module via $\rho$.

For a given ring $R$, let ${\bf Comp}(R)$ 
be the subcategory of ${\bf Comp}$ with objects $\{R, \mathbb{F}\}$ 
and morphisms all morphisms of the form $\{\operatorname{id}_R, \psi\}$ 
where $\psi$ is a map of chain complexes
over $R$. A ring homomorphism 
$\rho: R \rightarrow S$, 
naturally induces the 
\emph{base change functor}
\[
\rho_{*}\: {\bf Comp}(R) \rightarrow {\bf Comp}(S)
\]
given by 
\[
\rho_* 
\{R, \mathbb{F}  \} 
= \{S, S \otimes_{\rho} \mathbb{F}\} 
\]
and
\[
\rho_*
\{\operatorname{id}_R, \psi\} = 
\{\operatorname{id}_S, \operatorname{id}_S \otimes_\rho \psi\}.  
\]

\begin{dfn}
Let $\chi$ be a subcategory of 
{\bf Comp}. 

(a) 
We say that an object $\{R, \mathbb{F}\}$ of $\chi$ is an \emph{$R$-object of $\chi$}. 

(b)
We write $\chi(R)$ for the subcategory of 
{\bf Comp} with objects the $R$-objects of $\chi$ and morphisms those morphisms of $\chi$ that are also morphisms of {\bf Comp}($R$).

(c)
A  \emph{closed subcategory} is a 
subcategory $\chi$ of {\bf Comp} which is 
closed under base change in the sense that
for each ring homomorphism 
$\rho\: R\rightarrow S$ the base change 
functor 
$
\rho_*\: 
\text{\bf Comp}(R)\rightarrow
\text{\bf Comp}(S)
$
induces by restriction a 
functor (also called a base change functor) 
$
\rho_*\: \chi(R)\rightarrow \chi(S).   
$
\end{dfn}

\begin{definition}
Let $\chi$ be a closed subcategory, and  
let  
$\mathcal{F}: \chi \rightarrow \bf Comp$ 
be a functor. 

(a) 
We say that $\mathcal F$ is  
\emph{layered} if the following conditions 
hold:
\begin{itemize}
\item
For each $R$ the functor 
$\mathcal F$ restricts to  
a functor 
$
\mathcal F_R\: 
\chi(R) \longrightarrow 
\text{\bf Comp}(R).  
$
In particular, for each complex $\mathbb F$ 
of free $R$-modules we have 
$
\mathcal F\{R,\mathbb F\}= 
\{R, \mathcal F_R(\mathbb F)\},  
$
where $\mathcal F_R(\mathbb F)$ is again 
a chain complex of free $R$-modules. 
%

\item
For every ring homomorphism 
$\rho: R \rightarrow S$ there exists an 
isomorphism of functors \ 
$
\beta_{\rho}\: 
\mathcal{F}_S \circ \rho_{*} 
\longrightarrow 
{\rho_{*} \circ \mathcal{F}_R}, 
$ 
\ i.e. $\mathcal{F}$ commutes with base change. 
\end{itemize}

(b) 
We say that $\mathcal F$ 
is \emph{rigid} if 
it is layered and also satisfies 
\begin{itemize}
\item
For any $i\ge 0$ and any object 
$\{R,\mathbb F\}$ of $\chi$ one has 
$\HH_i\mathcal F_R(\mathbb F)=0$ if 
and only if 
$\HH_j\mathcal F_R(\mathbb F)=0$ for 
all $j\ge i$. 
\end{itemize}
\end{definition}

\begin{dfn}
Let $\mathbb{K}$ be a ring, $R$ a 
polynomial ring over $\mathbb{K}$ in 
finitely many variables, $\chi$ a closed 
subcategory, and $A$ an object of 
$\chi(R)$. 

(a) 
$A$ is said to be a \emph{$\chi$-generic
object over $\mathbb{K}$} if for any
commutative $\mathbb{K}$-algebra $S$ and
any $S$-object $B$ of $\chi$ there is a
homomorphism of $\mathbb{K}$-algebras
$\rho: R \rightarrow S$ such that
$\rho_{*}(A)$ and $B$ are isomorphic in
$\chi(S)$.

(b) 
If for every ring $\mathbb{K}$ there exists a $\chi$-generic object over $\mathbb{K}$, the closed subcategory $\chi$ is said to be \emph{sufficiently generic}.
\end{dfn} 

The Rigidity Criterion below is used to prove the rigidity of certain complexes which are built from symmetric powers, exterior powers, or Schur functors of a module. 

\begin{prop}[Rigidity Criterion ~\cite{T2007}] \label{rigiditycriterion}
Let $\chi$ be a sufficiently generic closed subcategory, let 
\[
\mathcal{F}: \chi \rightarrow \bf{Comp}
\]
be a layered functor, and assume that for every ring $\mathbb{K}$ there exists a $\chi$-generic over $\mathbb{K}$ object $\{ R, \mathbb{F}\}$ such that the complex $\mathcal{F}_R(\mathbb{F})$ is acyclic. Then the functor $\mathcal{F}$ is rigid. 
\end{prop}


The following example is essential to the proofs of our main results.

\begin{Example}\label{E:sufficiently-generic} 
(a) 
Let $t$ and $s$ be positive integers. Let $\mathcal{M}_{s,t}$ be the full subcategory of {\bf Comp} with objects $\{R, \mathbb{F} \}$ such that $\operatorname{rank} F_0 = s$, $\operatorname{rank} F_1 = t$, and $F_i = 0$ for $i \geq 2$. In other words, $\mathcal{M}_{s,t}$ is the category of homomorphisms from free modules of rank $t$ to free modules of rank~$s$.
%
%
Let $R = \mathbb{K}[x_{i,j} | 1 \leq i \leq s, 1 \leq j \leq t ]$ be the polynomial ring over $\mathbb{K}$ on the indicated set of variables. Let 
\[
\mathbb{G} = 0 \rightarrow R^t \xrightarrow{\bf{X}} R^s \rightarrow 0
\]
be the complex with $R^s$ in homological degree $0$ and such that the matrix for the map in the standard bases of $R^t$ and $R^s$ is $\bf{X}$ $= (x_{ij})$. 
Let $S$ be a $\mathbb{K}$-algebra. An $S$-object of $\mathcal{M}_{s,t}$ is of the form $$0 \rightarrow S^t \xrightarrow{\bf{Y}} S^s \rightarrow 0.$$ Define $\rho: R \rightarrow S$ by sending $x_{ij}$ to $y_{ij}$, where $\bf{Y}$ $= (y_{ij})$. The base change functor $\rho_{*}$ sends 
$$
0 \rightarrow R^t \xrightarrow{\bf{X}} 
R^s \rightarrow 0
$$ 
to 
$$
0 \rightarrow S^t \xrightarrow{\bf{Y}} 
S^s \rightarrow 0. 
$$ 
Thus, $\{R, \mathbb{G} \}$ is an 
$\mathcal{M}_{s,t}$-generic over $\mathbb{K}$ 
object and $\mathcal{M}_{s,t}$ is a sufficiently 
generic closed subcategory.

(b) 
Let $\widetilde{\mathcal M}_{s,t}$ be the 
subcategory 
of $\mathcal M_{s,t}$ with same objects, and 
morphisms 
those morphisms $\{\rho, \phi_\bullet\}$ 
such that 
$\phi_1$ is an isomorphism. In particular, 
isomorphisms 
in $\mathcal M_{s,t}$ are isomorphisms in 
$\widetilde{\mathcal M}_{s,t}$. 
It is straightforward 
to verify that the $\mathcal M_{s,t}$-generic 
over 
$\mathbb K$ object from part (a) is also an 
$\widetilde{\mathcal M}_{s,t}$-generic over 
$\mathbb K$ 
object, thus $\widetilde{\mathcal M}_{s,t}$ is 
also a sufficiently generic closed subcategory. 
\end{Example}

For the rest of this paper, 
the phrase ``generic case" 
refers to the case where 
\smash[t]{
$
\mathbb{F}: 
0 \rightarrow F \xrightarrow{\phi} G~\rightarrow~0
$
}
is a complex of free $R$-modules with 
$f=\rank F$ and $g=\rank G$, and the pair 
$\{R, \mathbb{F}\}$
is the $\mathcal{M}_{g,f}$-generic over 
$\mathbb{K}$ object from 
Example~\ref{E:sufficiently-generic}(a) for some 
$\mathbb{K}$.

\section{Schur and Weyl Modules}
\label{S:Schur-Weyl}

We review basic facts about Schur and Weyl 
modules. For more details and proofs the 
reader is referred to \cite{ABW1982} and to 
the excellent exposition in \cite{W2003}.

Recall that a \emph{partition} is a sequence 
$\lambda = (\lambda_1, \lambda_2, \dots)$ of 
nonnegative integers such that 
$\lambda_1 \geq \lambda_2 \geq \cdots$ and 
only a finite number of elements in the 
sequence are nonzero. The \emph{weight} of 
the partition $\lambda$ is the integer 
$|\lambda|=\sum_i\lambda_i$. 
We denote by  
$
\tilde{\lambda} = 
(\tilde{\lambda}_1, \tilde{\lambda}_2, \dots )
$ 
the \emph{conjugate} to $\lambda$ partition, 
where 
$\tilde{\lambda}_i$ is the number of 
$\lambda_j$'s that are $\geq i$.

A \emph{skew 
partition} $\lambda/\mu$ is a pair of 
partitions 
$\lambda = (\lambda_1, \lambda_2, \dots)$ and 
$\mu = (\mu_1, \mu_2, \dots)$ such 
that $\mu_i \leq \lambda_i$ for all $i$. In this 
case, we also write $\mu \subseteq \lambda$. 
We set $|\lambda/\mu|=\sum_i\lambda_i-\mu_i$. 
Let $(a_{ij})$ be the matrix 
defined as 
$$
a_{ij} =
\begin{cases}
1 & \text{if } \mu_i + 1 \leq j \leq \lambda_i; \\
0 & \text{otherwise}.
\end{cases}
$$
\begin{definition}
Let $F$ be a free $R$-module. The 
\emph{Schur module} $L_{\lambda/\mu}F$ 
is the 
image of the following composition 
$d_{\lambda/\mu}(F)$ of 
maps: 
\[
\begin{CD}
\wedge^{\lambda/\mu} =
\wedge^{\lambda/\mu}F := 
\wedge^{\lambda_1-\mu_1} F \otimes \wedge^{\lambda_2-\mu_2} F \otimes \cdots \\
@VV \Delta \otimes \Delta \otimes \cdots V \\
(\wedge^{a_{11}} F \otimes \wedge^{a_{12}} F \otimes \cdots) \otimes (\wedge^{a_{21}} F \otimes \wedge_{a_{22}} F \otimes \cdots) \otimes \cdots \\
@| \\
(S_{a_{11}} F \otimes S_{a_{12}} F \otimes \cdots) \otimes (S_{a_{21}} F \otimes S_{a_{22}} F \otimes \cdots)  \otimes \cdots \\
@VV \text{Rearrange Factors}V \\
(S_{a_{11}} F \otimes S_{a_{21}} F \otimes \cdots) \otimes (S_{a_{12}} F \otimes S_{a_{22}} F \otimes \cdots) \otimes \cdots \\
@VV \text{Multiplication}V \\
S_{\lambda/\mu} =
S_{\lambda/\mu}F :=
S_{\tilde{\lambda}_1-\tilde{\mu_1}} F \otimes S_{\tilde{\lambda}_2-\tilde{\mu_2}}  F \otimes \cdots . 
\end{CD}
\]
\end{definition}

\begin{Example}
Let $F$ be a free $R$-module, let 
$\lambda = (2,1)$, and take $\mu=0$. 
Then we have 
$\tilde{\lambda} = (2,1)$ and 
\[
(a_{ij})= \begin{bmatrix} 1 & 1 \\ 1 & 0 \end{bmatrix}.
\]
%
%
The map $d_{\lambda/\mu}(F)$ is the composition
\[
{\begin{CD}
\wedge^2 F \otimes F \\
@VV\Delta \otimes \Delta V \\
(F \otimes F) \otimes F \\
@VV \text{Rearranging Terms}V \\ 
(F \otimes F) \otimes F \\
@VV \text{Multiplication}V \\
S_2(F) \otimes F
\end{CD}}
\]
%
%
Given $a,b,c \in F$, under this composition 
we have 
\begin{align*}
 (a \wedge b) \otimes c  &\mapsto 
 (a\otimes b - b \otimes a) \otimes c 
= a \otimes b \otimes c - b \otimes a \otimes c  
\\ &\mapsto a \otimes c \otimes b - 
            b \otimes c \otimes a
\\ &\mapsto ac \otimes b - bc \otimes a. 
\end{align*}
\end{Example}

Note that when $\mu=0$ and $\lambda$
is a partition of the form 
$\lambda=(1,1,\dots,1)$ with $k$ number of $1$'s, 
the Schur module $L_{\lambda}F=L_{\lambda/0}F$ 
is equal to the $k$th symmetric power $S_kF$; 
while for a partition $\lambda$ 
of the form $\lambda=(k)$, the Schur module
$L_{\lambda}F$
is equal to the $k$th exterior power 
$\wedge^k F$.

\begin{definition}\label{D:skew-shape}
The \emph{skew-shape} $\Delta_{\lambda/\mu}$ 
associated with $\mu\subseteq \lambda$ is the set 
of integer pairs 
$
\Delta_{\lambda/\mu} = 
\{ (i,j) \mid  a_{ij}=1 \}. 
$
It is often visualized by replacing the 
nonzero entries in the matrix $(a_{ij})$ 
with rectangular boxes, and removing all zero 
entries. For example, for $\lambda=(4,3,2)$
and $\mu = (3,1)$, the skew shape looks
like this: 
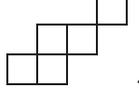
\[
\young(:::~,:~~,~~) \ 
.
\]
\end{definition}

\begin{dfn}
Let $S$ be a 
set, let $\lambda$ and $\mu$ be partitions 
with $\mu \subseteq \lambda$, and let 
$\Delta_{\lambda/\mu}$ be the skew-shape 
associated to this pair. A 
\emph{tableau} of shape
$\lambda/\mu$ with values in the set $S$ is
a function from $\Delta_{\lambda/\mu}$ to
$S$. The set of all such tableaux is
denoted by 
$\operatorname{Tab}_{\lambda/\mu}(S)$.
\end{dfn}

\begin{example}
Let $\lambda = (4,3,2)$, $\mu = (3,1)$, and 
$S = \{a ,b, c, d\}$. 
Then $\Delta_{\lambda/\mu}$ looks like
$$
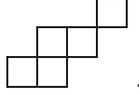
\young(:::~,:~~,~~) \ 
.
$$
In terms of coordinates we have 
$
\Delta_{\lambda/\mu} = 
\{(1,4), (2,2), (2,3), (3,1), (3,2)\}.  
$
Now define a map to $S$ by
\begin{align*}
(1,4) &\mapsto a\\
(2,2) &\mapsto b \\
(2,3) &\mapsto c \\
(3,1) &\mapsto c \\
(3,2) &\mapsto d
\end{align*}
This function is a tableau of shape 
$\lambda/\mu$ with values in $S$. We may think 
of it as a way to fill the tableau with values 
in $S$ as seen below
$$
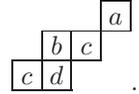
\young(:::a,:bc,cd).
$$
\end{example} 

\begin{definition}
Let $F$ be a free $R$-module,  
and let 
$S$ be a subset of $F$. 
Given an element 
$
T \in 
\operatorname{Tab}_{\lambda/\mu}S,
$ 
we associate to it a simple tensor  
$Z_T \in \wedge^{\lambda/\mu}F$, 
whose component  in 
$\wedge^{\lambda_j-\mu_j}F$ is the exterior 
product of the elements in the $j$th row of $T$. 
For example, the tableau above corresponds 
to the simple tensor 
$a \otimes (b \wedge c) \otimes (c \wedge d)$ in  
$
\wedge^{\lambda/\mu} = 
\wedge^1 \otimes \wedge^2 \otimes \wedge^2. 
$
\end{definition}

The Weyl module $K_{\lambda/\mu}F$ of a free 
module $F$ is defined in a dual manner to that of 
the Schur 
module, and is the image of a map 
$
d'_{\lambda/\mu}: 
D_{\lambda/\mu}F \rightarrow 
\wedge_{\tilde{\lambda}/\tilde{\mu}}F.
$ 
Explicitly, the Weyl module is defined 
as the image of the following composition
of maps:
%
\[
\begin{CD}
D_{\lambda_1-\mu_1} F \otimes D_{\lambda_2-\mu_2} F \otimes \cdots \\
@VV \Delta \otimes \Delta \otimes \cdots V \\
(D_{a_{11}} F \otimes D_{a_{12}} F \otimes \cdots) \otimes (D_{a_{21}} F \otimes D_{a_{22}} F \otimes \cdots) \otimes \cdots \\
@| \\
(\wedge^{a_{11}} F \otimes \wedge^{a_{12}} F \otimes \cdots) \otimes (\wedge^{a_{21}} F \otimes \wedge^{a_{22}} F \otimes \cdots)  \otimes \cdots \\
@VV \text{Rearrange Factors}V \\
(\wedge^{a_{11}} F \otimes \wedge^{a_{21}} F \otimes \cdots) \otimes (\wedge^{a_{12}} F \otimes \wedge^{a_{22}} F \otimes \cdots) \otimes \cdots \\
@VV \text{Multiplication}V \\
\wedge^{\tilde{\lambda}_1-\tilde{\mu_1}} F \otimes \wedge^{\tilde{\lambda}_2-\tilde{\mu_2}}  F \otimes \cdots . 
\end{CD}
\]
Note that the Schur module of a free module $F$ is
zero if and only if a row in the skew partition 
$\lambda/\mu$ is longer 
than the rank of $F$, whereas the Weyl module is 
zero if and only if a column is longer than the rank.

To conclude this section, 
we recall how to extend the definition of
the Schur module $L_{\lambda/\mu}M$ to 
the case when 
$M$ which is not necessarily free.

Let $p_i = \lambda_i - \mu_i$, let 
$q=\max\{i \mid p_i\ne 0\}$, and define the functor $\wedge^{(\lambda/\mu)+}$ to be the direct sum
\[
\sum_{i=1}^{q-1} \ \ \sum_{t = \mu_i - \mu_{i+1}}^{\lambda_{i+1}-\mu_{i+1}} 
\wedge^{p_1} \otimes \cdots \otimes \wedge^{p_{i-1}} \otimes \wedge^{p_i +t} \otimes \wedge^{p_{i+1} -t} \otimes \wedge^{p_{i+2}} \otimes \cdots 
\wedge^{p_q} 
\]
Define a map 
\[
\square_{\lambda/\mu}M:
\wedge^{(\lambda/\mu)+} M \rightarrow 
\wedge^{\lambda/\mu} M
\]
as follows. 
If $q < 2$, set 
$\wedge^{(\lambda/\mu)+} = 0$. If $q = 2$, 
the map 
\[
\square_{\lambda/\mu}M : 
\sum_{t=\mu_1 - \mu_2+1}^{\lambda_2 - \mu_2} \wedge^{p_1 + t} M \otimes \wedge^{p_2 - t}M \rightarrow \wedge^{p_1} M \otimes \wedge^{p_2} M
\]
is defined on the component 
$\wedge^{p_1+t}M\otimes\wedge^{p_2-t}M$
as the composition 
\[
\wedge^{p_1+t}M\otimes\wedge^{p_2-t}M 
\xrightarrow{\Delta\otimes 1} 
\wedge^{p_1}M\otimes\wedge^{t}M
             \otimes\wedge^{p_2-t}M
\xrightarrow{1\otimes\wedge}
\wedge^{p_1}M\otimes\wedge^{p_2}M
\]
where $\Delta$ and $\wedge$ are the diagonal 
and exterior multiplication maps, respectively. 
If $q>2$,  set 
$\lambda^i = (\lambda_i, \lambda_{i+1})$, $\mu^i = (\mu_i, \mu_{i+1})$ and define $\square_{\lambda/\mu}M$, to be 
$$
\sum_{i=1}^{q-1} 
1_1 \otimes \cdots \otimes 1_{i-1} \otimes \square_{\lambda^i/\mu^i}M \otimes 1_{i+2} \otimes \cdots \otimes 1_q.
$$
Now we have  

\begin{definition}[\cite{ABW1982}] 
Let $M$ be an 
$R$-module, and let 
$\mu\subseteq\lambda$ be partitions. 
The $R$-module 
$
L_{\lambda/\mu}M = 
\Coker(\square_{\lambda/\mu})
$ 
is called the \emph{Schur module} of $M$. 
\end{definition}

\section{Schur Complexes}
\label{prelimsymmetric}

In this section 
we recall some basic properties of Schur 
complexes. All definitions and results 
presented here are 
due to \cite{ABW1982}, to which we 
refer the reader for more 
details and proofs. 

Let $M$ be a module over $R$ and let $k\geq0$.
Let 
$$
\phi: F \rightarrow G
$$ 
be a homomorphism of finite free $R$-modules 
with cokernel $M$, 
let $f=\rank F$, and let $g=\rank G$. From the 
definition, it is straightforward that the maps 
below have cokernels $S_k M$ and $\wedge^k M$ 
respectively: 
\begin{equation}
\label{E:symmetric} 
F \otimes S_{k-1}G \longrightarrow S_k G,   
\quad  
f \otimes u \longmapsto \phi(f)u,  
\end{equation}
\begin{equation} 
\label{E:exterior}
F \otimes \wedge^{k-1} G \longrightarrow 
          \wedge^k G, 
\quad  
f \otimes v \longmapsto \phi(f)\wedge v,
\end{equation}
where $f \in F, u \in S_{k-1} G$, and 
$v \in \wedge^{k-1} G$. 

%
%

For each $t$, we have a complex $S_t \phi$, 
referred to by \cite{E1995} as the $t$th graded 
strand of the Koszul complex of $\phi$ 
\[ 
0 \rightarrow 
\wedge^n F \otimes S_{t-n}G 
\xrightarrow{d_n} 
\wedge^{n-1}F \otimes S_{t-n+1}G 
\xrightarrow{d_{n-1}} \dots \xrightarrow{d_2} 
F \otimes S_{t-1}G 
\xrightarrow{d_1} S_t G \rightarrow 0 
\]
with
$n= \text{min}\{\text{rank} F, t\}$
where 
\[
d_j(
f_1 \wedge f_2 \wedge \dots \wedge f_j \otimes m
) 
= 
\sum^j_{r=1} 
(-1)^{r-1}f_1 \wedge \dots \wedge 
    \hat{f_r} \wedge \dots \wedge 
f_j \otimes \phi(f_r)m.
\]
We have by \eqref{E:symmetric} that 
$H_0 (S_t \phi) = S_t M $. 
%
%
%
%
%
%
Similarly, Lebelt~\cite{L1975} studied a complex 
$$ 
\wedge^t \phi :\quad  
0 \rightarrow D_t F 
  \rightarrow D_{t-1}F \otimes G 
  \rightarrow \cdots \rightarrow \wedge^t G \rightarrow 
0
$$
where 
for $0 \leq i \leq t$ the free module  
$D_i F = S_i(F^*)^*$  is the 
$i$th divided power of $F$, 
see \cite[Section A2.4]{E1995}, 
and
the  differential is given by 
\[
b_1^{(a_1)}b_2^{(a_2)}\dots\ b_f^{(a_f)} 
\otimes v 
\longmapsto 
\sum^f_{r=1} 
b_1^{(a_1)} \dots\  
b_r^{(a_r-1)} \dots\   
b_f^{(a_f)} \otimes
\bigl(\phi(b_r)\wedge v\bigr).
\]
By \eqref{E:exterior}, this complex has 
$\wedge^t M$ as it zeroth homology.
%
%

Set
\[
\wedge^{\lambda/\mu} \phi = 
\wedge^{\lambda_1-\mu_1} \phi \otimes 
\wedge^{\lambda_2 - \mu_2}\phi \otimes 
\cdots 
\]
and 
\[
S_{\tilde{\lambda}/\tilde{\mu}} \phi= 
S_{\tilde{\lambda_1}-\tilde{\mu}} \phi \otimes 
S_{\tilde{\lambda_2} - \tilde{\mu_2}} \phi\otimes 
\cdots . 
\] 

\begin{dfn}
The \emph{Schur complex} $L_{\lambda/\mu}\phi$ is the image of the map 
\[
d_{\lambda/\mu} \phi: 
\wedge^{\lambda/\mu} \phi 
\longrightarrow S_{\tilde\lambda/\tilde\mu} \phi
\] 
defined as the composition
\[
\begin{CD}
\wedge^{\lambda_1-\mu_1} \phi \otimes \wedge^{\lambda_2-\mu_2} \phi \otimes \cdots \\
@VV \Delta \otimes \Delta \otimes \cdots V \\
(\wedge^{a_{11}} \phi \otimes \wedge^{a_{12}} \phi \otimes \cdots) \otimes (\wedge^{a_{21}} \phi \otimes \wedge^{a_{22}} \phi \otimes \cdots) \otimes \cdots \\
@| \\
(S_{a_{11}} \phi \otimes S_{a_{12}} \phi \otimes \cdots) \otimes (S_{a_{21}} \phi \otimes S_{a_{22}} \phi \otimes \cdots)  \otimes \cdots \\
@VV \text{Rearrange Factors}V \\
(S_{a_{11}} \phi \otimes S_{a_{21}} \phi \otimes \cdots) \otimes (S_{a_{12}} \phi \otimes S_{a_{22}} \phi \otimes \cdots) \otimes \cdots \\
@VV \text{Multiplication}V \\
S_{\tilde{\lambda}_1-\tilde{\mu_1}} \phi \otimes S_{\tilde{\lambda}_2-\tilde{\mu_2}}  \phi \otimes \cdots 
\end{CD}
\]
where $\Delta$ is the diagonal map.  
\end{dfn}


\begin{remark} \label{layered}
From the definition of the complex  
$L_{\lambda/\mu}\phi$ it is straightforward 
that it 
commutes with base change,  hence  
induces a layered functor 
$
L_{\lambda/\mu}\: 
\mathcal M_{g,f} \rightarrow \text{\bf Comp}.
$
\end{remark}


\begin{dfn}
Let $S$ be a totally ordered set, let $X$ be a subset of $S$, and let 
$T\in\operatorname{Tab}_{\lambda/\mu}(S)$. 

(a) 
The  tableau $T$ is said to be 
$\emph{row-standard}$ mod $X$ if each row of $T$ 
is non-decreasing, and if, when repeats 
occur in a row, they occur only among 
elements of $X$. Call 
$T$  $\emph{column-standard}$ mod $X$ if 
each column is non-decreasing, and if, when 
repeats occur in a column, they occur only 
among 
elements in the complement of $X$. Call $T$  
$\emph{standard}$ mod $X$ if T is both 
row- and column-standard mod $X$. 

(b) 
Let $\phi\: F \rightarrow G$ be a map of finite 
free $R$-modules, and suppose that $X$ is a 
subset 
of $F$ and $S\setminus X$ is a subset of $G$. 
When $T$ is standard mod $X$, write 
$Z_T$ for the simple tensor in 
$\wedge^{\lambda/\mu}\phi$ 
whose component in 
$\wedge ^{\lambda_i-\mu_i} \phi$ is given by 
the product in the divided powers algebra of 
$F$ of the elements of $X$ in row $j$ (a 
repeated element is given its corresponding 
divided power) tensored by 
the exterior product of the elements from
$S\setminus X$ that are in row $j$ of $T$; 
see the examples below.
\end{dfn}

\begin{Example}
Let $F$ and $G$ be free modules with bases $X = \{a, b\}$ and $Y = \{x, y, x\}$, respectively. Let $S = \{a, b, x, y, z\}$ be ordered by  
$a<b< x<y<z$, let $\lambda = (4,3)$ and let 
$\mu = (2,1)$. 
The following are examples of standard tableaux 
$T$ mod $X$ and their corresponding simple  
tensors $Z_T$ in  
$\wedge^{\lambda/\mu}\phi = 
\wedge^2\phi\otimes\wedge^2\phi$: 
\begin{align*}
T   &=\young(::aa,:ab) \\ 
Z_T &= (a^{(2)}\otimes 1)\otimes (ab\otimes 1) 
       \in  
       (D_2F\otimes\wedge^0G)\otimes
       (D_2F\otimes \wedge^0G)   
       \subset  
       \wedge^2\phi\otimes\wedge^2\phi; 
\end{align*}
\begin{align*}
T   &= \young(::xy,:ax) \\ 
Z_T &= (1\otimes x \wedge y)\otimes 
       (a \otimes x) 
       \in  
       (D_0F \otimes \wedge^2 G) \otimes 
       (D_1F \otimes \wedge^1 G) 
       \subset 
       \wedge^2 \phi \otimes \wedge^2 \phi; 
\end{align*}
\begin{align*} 
T   &= \young(::ab,:yz) \\ 
Z_T &= (ab\otimes 1) \otimes (1\otimes y\wedge z)
       \in 
       (D_2F \otimes \wedge^0 G) \otimes 
       (D_0F \otimes \wedge^2G) 
       \subset 
       \wedge^2 \phi \otimes \wedge^2 \phi. 
\end{align*} 
\end{Example}

The following theorem provides a free basis for a Schur complex $L_{\lambda/\mu} \phi$.

\begin{thm}\label{T:schur-basis}
Let $\lambda$ and $\mu$ 
be partitions with $\mu \subset \lambda$, and 
let $\phi: F \rightarrow G$ be a map of free 
$R$-modules. Let $X$ and $Y$ be bases for $F$ and 
$G$ respectively, 
and let $S = X \sqcup Y$ be totally ordered. 
Then the set 
\[
\{ 
d_{\lambda/\mu}(Z_T) \mid  
T \in \operatorname{Tab}_{\lambda/\mu}(S) 
\text{ is standard mod } X 
\}
\] 
is a basis for $L_{\lambda/\mu}\phi$.
\end{thm}


We will also need the following results. 

\begin{thm} \label{mapsplits}
If $\phi = \phi_1 \oplus \phi_2$, where $\phi_1$ is an isomorphism, then the natural inclusion map  
$
L_{\lambda/\mu}\phi_2\lra L_{\lambda/\mu}\phi
$ 
splits and yields an isomorphism 
$
L_{\lambda/\mu}\phi\cong L_{\lambda/\mu}\phi_2 \oplus E
$ 
where $E$ is a contractible
chain complex.  
\end{thm}

\begin{thm}
Let $\phi:F \rightarrow G$ be a map of free 
$R$-modules and let $\mu \subseteq \lambda$ be 
partitions. Define $(L_{\lambda/\mu} \phi)_j$ 
to be the component in homological 
degree $j$ of the complex 
$L_{\lambda/\mu}\phi$. There is a natural 
filtration on $(L_{\lambda/\mu}\phi)_j$ whose 
associated graded module is 
$$
\sum_{\mu \subseteq \gamma \subseteq \lambda;            |\lambda|-|\gamma|=j} 
L_{\gamma/\mu}G \otimes K_{\lambda/\gamma}F. 
$$
In particular, 
$(L_{\lambda/\mu}\phi)_0=L_{\lambda/\mu}G$. 
\end{thm}

\begin{corollary}\label{C:schur_complex_bounds} 
Let $\nu''$ be the partition with 
$\nu''_i=\min\{\lambda_i, \mu_i+\rank G\}$, and let 
$\nu'$ be the partition with 
$
\tilde\nu'_i=
\max\{\tilde\mu_i, \tilde\lambda_i-\rank F\}. 
$

(a) 
$(L_{\lambda/\mu}\phi)_j\ne 0$ if and only if for some 
$\gamma$ with 
$\mu\subseteq\gamma\subseteq\lambda$ and
$|\lambda|-|\gamma|=j$ we have 
$\tilde\lambda_i -\tilde\gamma_i\le \rank F$ for all $i$, and 
$\gamma_t-\mu_t\le \rank G$ for all $t$, or equivalently,  
$\nu'\subseteq\gamma\subseteq\nu''$. 

(b) In particular, 
$L_{\lambda/\mu}\phi\ne0$ if and only if 
$\nu'\subseteq\nu''$. 
In that case:  

(c)
The unique maximal $j$ with 
$(L_{\lambda/\mu}\phi)_j\ne 0$ 
is obtained when taking $\gamma=\nu'$, hence equals 
\[
|\lambda| - |\nu'| = 
|\lambda| - 
\sum_i\max\{\tilde\mu_i,
            \tilde\lambda_i-\rank F\}. 
\]
The unique minimal $j$ with 
$(L_{\lambda/\mu}\phi)_j\ne 0$ 
is obtained when taking $\gamma=\nu''$, hence equals 
\[
|\lambda| - |\nu''| = 
|\lambda| - 
\sum_i\min\{\lambda_i, \mu_i + \rank G\}. 
\]

(d)
$(L_{\lambda/\mu}\phi)_j\ne 0$ if and only if  
\[
|\lambda|-\sum_t\min\{\lambda_t,\mu_t+\rank G\}
\le j\le 
|\lambda| - \sum_i\max\{\tilde\mu_i, \tilde\lambda_i-\rank F\}.
\] 
\end{corollary}

\begin{prop}\label{P:Schur-homology}
Let $\pi\: G \longrightarrow M=G/\im\phi$ be the 
canonical projection. Then 
\[
(L_{\lambda/\mu}\phi)_1 
\xrightarrow{\ d_1\ }  
L_{\lambda/\mu}G 
\xrightarrow{L_{\lambda/\mu}(\pi)}
L_{\lambda/\mu}M 
\longrightarrow 0 
\]
is exact, where $d_1$ is the differential 
of $L_{\lambda/\mu}\phi$. In particular, 
$H_0(L_{\lambda/\mu} \phi)= L_{\lambda/\mu}M$. 
%
\end{prop}

\section{More on Schur Complexes}
\label{S:threshhold}

The following notation and terminology 
will be useful for us  in the sequel. 
Our main point here is the introduction 
of the \emph{threshold number} in 
Definition~\ref{D:threshhold}. 

\begin{definition}
(a) 
We call the integer 
$
W(\lambda/\mu)=
\max\{\lambda_i-\mu_i\mid i\ge 1\}
$ 
the \emph{width} of $\lambda/\mu$, and we call 
$
H(\lambda/\mu)=
\max\{\tilde\lambda_i-\tilde\mu_i\mid 
       i\ge 1\}
$ 
the \emph{height} of $\lambda/\mu$. 

(b) 
For each integer $n$ define the partition 
$\nu''=\nu''(\lambda, \mu, n)$ by setting 
\[
\nu''_i=\min\{\lambda_i, \max\{\mu_i, \mu_i + n\}\}, 
\] 
and define the partition $\nu'=\nu'(\lambda,\mu,n)$ by setting 
\[
\tilde\nu_i'=\max\{\tilde\mu_i, \min{\{\tilde\lambda_i, \tilde\lambda_i - n}\}\}.
\]
\end{definition}

\begin{remark}
It is immediate from the definition 
that we have the following inclusions 
and equalities: 
{\small
\begin{align*}
\lambda = \dots = 
\nu'(\lambda,\mu,-1) = 
\nu'(\lambda,\mu, 0) &\ge \dots \ge  
\nu'(\lambda,\mu, H) = 
\nu'(\lambda,\mu, H + 1) = \dots = \mu; \\
\mu = \dots = 
\nu''(\lambda,\mu,-1) = 
\nu''(\lambda,\mu, 0) &\le \dots \le 
\nu''(\lambda,\mu, W) = 
\nu''(\lambda,\mu, W + 1) = \dots = 
\lambda; 
\end{align*}
}
\negthickspace 
where $H=H(\lambda/\mu)$ and 
$W=W(\lambda/\mu)$. Furthermore, 
if $\mu\ne\lambda$ then all inequalities 
above are strict. 
\end{remark}

\begin{definition}
\label{D:threshhold} 
(a) 
For any integers $f$ and $g$ we 
set 
\[
T_{\lambda/\mu}(f,g) = 
\max\{
t\mid \nu'(\lambda,\mu,f-t)\subseteq 
\nu''(\lambda,\nu,g-t)\},  
\]
and call this the \emph{threshold number} 
of $\lambda/\mu$ with respect to the pair 
$(f,g)$. 

(b) 
For any integer $n$ we set 
\begin{align*}
l_n &= |\nu'(\lambda,\mu,n - 1)| - 
    |\nu'(\lambda,\mu,n)|
\quad\text{and}\quad \\ 
k_n &=|\nu''(\lambda,\mu,n)| - 
    |\nu''(\lambda,\mu, n-1)|. 
\end{align*}
\end{definition}

As a straightforward consequence of the 
definitions we have the following basic 
properties: 

\begin{remarks}\label{R:threshold}
Let $H=H(\lambda/\mu)$, let  
$W=W(\lambda/\mu)$, and let 
$T=T_{\lambda/\mu}(f,g)$. 

(a) 
$T_{\lambda/\mu}(f,g)\ge 0$ if and only if 
$
\nu'(\lambda,\mu,f)\subseteq
\nu''(\lambda,\mu,g).
$

(b) 
$T_{\lambda/\mu}(f-1,g-1)=T_{\lambda/\mu}(f,g)-1$. 

(c) 
$W\ge g - T$. Indeed, by the definition of 
$T$ we have 
$\nu'(\lambda,\mu,f-T)\subseteq \nu''(\lambda,\mu,g-T)$, 
and 
$
\nu'(\lambda,\mu,f-T-1)\not\subseteq
\nu''(\lambda,\mu,g-T-1). 
$
If $W<g-T$ then 
$
\nu''(\lambda,\mu,g-T-1)=
\nu''(\lambda,\mu,g-T) = \lambda
$
and so 
$
\nu'(\lambda,\mu,f-T-1)\subseteq 
\nu''(\lambda,\mu,g-T-1), 
$
a contradiction. 

(d)
Either $T_{\lambda/\mu}(0,g)<0$ or 
$T_{\lambda/\mu}(0,g)=g-W\ge 0$. 

(e)
Either $T_{\lambda/\mu}(f,0)<0$ or 
$T_{\lambda/\mu}(f,0)=f- H\ge 0$. 

(f)
If $n\ge 1$ then 
$k_n$ is the number of $t$'s such that 
$\mu_t+ n \le \lambda_t$, and  
$l_n$ is the number of $t$'s such that 
$\tilde\mu_t\le\tilde\lambda_t-n$. In particular, 
$\lambda$ differs from $\mu$ in exactly 
$k_1$ of its rows and in exactly $l_1$ 
of its columns. 

(g)   
$k_n=0$ for $n\le 0$ and $n\ge W$, and 
$l_n=0$ for $n\le 0$ and $n\ge H$. 

(h)
We have $l_1\ge \dots \ge l_{H-1}$ and 
$k_1\ge \dots \ge k_{W-1}$. 

(i) 
For each $n$ we have 
\[
\sum_{t\ge n+1} k_t = |\lambda|-|\nu''(\lambda,\mu,n)|
\quad\text{and}\quad 
\sum_{t\le n} l_t = |\lambda| - |\nu'(\lambda,\mu,n)|. 
\]
In particular, 
\[
|\lambda|-|\mu| = 
\sum_{n=1}^W k_n = 
\sum_{n=1}^H l_n.  
\]

\end{remarks}

We can now reformulate Corollary~\ref{C:schur_complex_bounds} as follows. 

\begin{corollary} \label{nonzeroiff}
The complex $L_{\lambda/\mu}\phi$ is nonzero if and only 
if 
\[
\nu'(\lambda,\mu,\rank F)\subseteq  
\nu''(\lambda,\mu,\rank G), 
\]
which is if and only if \ 
$T_{\lambda/\mu}(\rank F,\rank G)\ge 0$. 
In that case the component $(L_{\lambda/\mu}\phi)_j$ is nonzero if and only if 
\[
\sum_{t\ge 1+\rank G} k_t \ \le \ j \ \le 
\sum_{1\le t\le\rank F} l_t. 
\]
\end{corollary}

We will need the following related basic
observation.

\begin{proposition} \label{P:nonzero-differential}
Suppose $\mu\subsetneq\lambda$ and let 
$\phi\: F\lra G$ be a generic map. The   
$j$-th differential of 
$L_{\lambda/\mu}\phi$ 
satisfies $\partial_j\ne 0$  
whenever \ 
$
\start L_{\lambda/\mu}\phi < j \le 
\finish L_{\lambda/\mu}\phi.
$
\end{proposition}

\begin{proof}
Let $f= \rank F$ and $g=\rank G$, let 
$X=\{x_1,\dots, x_f\}$ be a basis for $F$, and 
let $Y=\{y_1,\dots,y_g\}$ be a basis for $G$. 
We order the disjoint union  $S=X\sqcup Y$ by 
setting $y_1<\dots< y_g < x_1 <\dots< x_f$. 

If $L_{\lambda/\mu}\phi=0$ the proposition is 
trivially true. Suppose 
$(L_{\lambda/\mu}\phi)_j\ne 0$. 
Thus there is a partition $\gamma$ such that 
$|\lambda|-|\gamma|=j$ and 
\[
\nu'=\nu'(\lambda,\mu,f)\subseteq \gamma \subseteq
\nu''(\lambda,\mu,g)=\nu''. 
\]
It suffices to show that if 
$j\ne\start L_{\lambda/\mu}\phi$ then we can  
specialize $\phi$ to a map $\psi\: F\lra G$ 
so that for the $j$th  
differential $\delta_j$ of $L_{\lambda/\mu}\psi$
we have $\delta_j\ne 0$.  
So, assume $\gamma<\nu''$, and take $\psi$ to be 
the map that sends $x_1$ to $y_g$, and all the 
other $x_i$ to $0$. 
Consider the basis element 
$d_{\lambda/\mu}(Z_\Gamma)$ of 
$L_{\lambda/\mu}(\psi)_j$, where   
$\Gamma$ is the standard mod $X$ tableau of shape 
$\lambda/\mu$ given by 
\[
\Gamma(i,j) = 
\begin{cases}
x_{i-\tilde\gamma_j} &\text{ if } 
                      (i,j)\in\lambda/\gamma; \\ 
y_{j-\mu_i}          &\text{ if } 
                      (i,j)\in\gamma/\mu,  
\end{cases}
\]
and $Z_\Gamma$ is the corresponding basis element 
of $\wedge^{\lambda/\mu}\psi$. Let 
$J=\{ \tilde\gamma_j + 1\mid j\ge 1\}$.  
Note that $1\in J$, and that $2\le i\in J$ exactly 
if $\gamma_i <\gamma_{i-1}$. Let 
\[
A=
\{i\in J\mid \gamma_i<\min(\lambda_i,\mu_i + g)\}. 
\] 
Note that, as $\gamma<\nu''$, the set $A$ is nonempty.  
Therefore,   
for each $i\in A$, by adding one 
box to the $i$th row of $\gamma$ we obtain a new  
partition $\gamma^{(i)}$. As 
$\gamma<\gamma^{(i)}\le\nu''$, we obtain a  
corresponding standard 
mod $X$ tableau $\Gamma^{(i)}$ 
of shape $\lambda/\mu$ by setting 
\[
\Gamma^{(i)}(p,q) = 
\begin{cases}
y_{g} &\text{ if } 
       (p,q)=(i,\gamma_i +1); \\ 
\Gamma(p,q) &\text{ otherwise}.  
\end{cases}
\]
Now it is straightforward 
to check that 
the differential 
$d_j$ of $\wedge^{\lambda/\mu}\psi$ sends the 
basis element $Z_\Gamma$ to 
\[
d_j(Z_\Gamma) = 
\sum_{i\in A} (-1)^{\gamma_i} Z_{\Gamma^{(i)}}, 
\]
hence 
$
\delta_j\bigl(d_{\lambda/\mu}(Z_\Gamma)\bigr)=
d_{\lambda/\mu}\bigl(d_j(Z_\Gamma)\bigr)=
\sum_{i\in A} 
(-1)^{\gamma_i} d_{\lambda/\mu}(Z_{\Gamma^{(i)}}) 
\ne 0
$
as the elements $d_{\lambda/\mu}(Z_{\Gamma^{(i)}})$ 
are part of a basis for $L_{\lambda/\mu}\psi$. 
Thus $\delta_j\ne 0$ as desired. 
\end{proof}

\section{Acyclicity of Schur complexes}
\label{prob.ch}

Here is the first main result of the 
paper. 

\begin{thm} \label{sqrt}
Let $\mu\subsetneq\lambda$, 
let $\phi: F\lra G$ be a 
map of finite free $R$-modules, 
let $f=\rank F$, \ 
let \ $g=\rank G$, \ and   
let \ $T=T_{\lambda/\mu}(f, g)$. 
\begin{enumerate} 
\item 
Suppose that \ 
$f-g< H(\lambda/\mu)$. \ 
Then \ $r_n\ge 0$ \ for each \ $n\ge 1$, \ and: 

(a) 
For \ 
$
\sum_{i\ge 1+g-T} k_i \ < \ n \ \le \ 
\sum_{i\le f-T} l_i  
$
\ we have 
\[
\sqrt{\vphantom{I_j(\phi)}I(d_n)}=
\sqrt{\vphantom{I_j(\phi)}I_{1+T}(\phi)}; 
\]

(b)
For \ 
$
\sum_{i\le j-1} l_i \ < \ n \ \leq \ 
\sum_{i\le j}l_i
$ 
\ with \ $j\ge f-T+1$ \ we have
\[
\sqrt{\vphantom{I_j(\phi)}I(d_n)} = 
\sqrt{I_{f-j+1}(\phi)}; 
\]

(c) 
For \ 
$
\sum_{i\ge 2+g-j} k_i \ < \ n \ \le \  
\sum_{i\ge 1+g-j} k_i 
$ 
\ with \ $j\le T$ \ we have 
\[
\sqrt{\vphantom{I_j(phi)}I(d_n)}=\sqrt{I_j(\phi)}; 
\]

\item 
Suppose that \ $f-g\ge H(\lambda/\mu)$. \ Then \ 
$\finish L_{\lambda/\mu}\phi=|\lambda/\mu|$, 
\ and: 

(a) 
For \ $1\le n\le |\lambda/\mu|$ \ with \  
$|\lambda/\mu|-n$ \ even we have \ $r_n>0$ \ and 
\[
I(d_n)=0; 
\]

(b) 
For \ $1\le n\le |\lambda/\mu|$ \ with \  
$|\lambda/\mu|-n$ \ odd we have 
\[
\sqrt{\vphantom{I_g(\phi)}I(d_n)}\supseteq 
\sqrt{I_g(\phi)}; 
\]
\end{enumerate}
where $r_n$ are the expected ranks of the differentials $d_n$ in $L_{\lambda/\mu}(\phi)$.
\end{thm}

\begin{proof}
We proceed by induction on 
$m=\min(f,g)$. When $m=0$ we have 
either $f=0$ or $g=0$. In the former case 
$L_{\lambda/\mu}\phi$ is concentrated in 
homological degree $0$ and by 
Remark~\ref{R:threshold}(d) 
either $T<0$ or $T=g-W\ge 0$, thus the theorem is 
trivially true. In  the latter case, 
$L_{\lambda/\mu}\phi$ is concentrated in 
homological degree $|\lambda/\mu|$ and by 
Remark~\ref{R:threshold}(e) 
either $T<0$ or $T=f-H\ge 0$, thus,  in 
view of Remark~\ref{R:threshold}(a) and 
Corollary~\ref{nonzeroiff}, the theorem is 
again trivially true.

Now assume $m\ge 1$.  
By a standard argument, 
we may assume that the ring $R$ is a polynomial ring
over $\mathbb{Z}$ with indeterminates the entries of
$\phi$.
Note that, by Theorem~\ref{mapsplits}, 
after inverting any 
$m$ by $m$ minor 
of $\phi$ we obtain 
\begin{equation}\label{E:splits}
L_{\lambda/\mu}\phi\cong
L_{\lambda/\mu}\phi'\oplus E
\end{equation}
where $E$ 
is a contractible complex of free modules, hence 
split exact, and $\phi': F'\rightarrow G'$ 
is a map of free modules with 
$f'=\rank F'=f-m$ and $g'=\rank G'=g-m$. 
Since $f'-g'=f-g$ and the expected rank of $d_n$ 
is at least the same as the expected rank 
$r'_n$ for the differential
$d'_n$ of $L_{\lambda/\mu}\phi'$, all expected 
rank inequalities claimed in  the theorem follow
from our induction hypothesis. 
Since $R$ is a 
domain and $I(d'_n)$ is the localization of
$I(d_n)$, our induction hypothesis also yields 
(2) part (a). 
When in case (2) we have $g=m$,  
and $L_{\lambda/\mu}\phi'$ is nonzero and 
concentrated in homological degree $|\lambda/\mu|$, in particular 
$\finish L_{\lambda/\mu}\phi=|\lambda/\mu|$.
Thus in (2) part (b) we have $r'_n<0$ and so 
$I(d'_n)$ equals the whole ring; therefore 
the radical of $I(d_n)$ contains $I_g(\phi)$, 
completing the proof of (2) part (b). 

It remains to prove the radical equalities 
claimed in case (1). We  
note that they are trivially satisfied 
when $T<0$ because then $L_{\lambda/\mu}\phi=0$ 
and so $I(d_n)=R$ for all $n$, and we assume for 
the rest of this proof that $T\ge 0$ and 
that we are in case (1). Thus, 
we already know that the expected rank $r_n$  
for the differential $d_n$ of 
$L_{\lambda/\mu}\phi$ is nonnegative for each  
$n\ge 1$. Furthermore, 
$L_{\lambda/\mu}\phi'$ is concentrated in 
homological degree $0$, hence 
$L_{\lambda/\mu}\phi$ becomes acyclic after 
localization and thus for each $n\ge 1$ the 
expected rank $r_n$ is the same as the 
determinantal rank of $d_n$. Therefore from 
Proposition~\ref{P:nonzero-differential} we obtain 
that $r_n=0$ for 
$1\le n\le\start L_{\lambda/\mu}\phi$ and 
$n> \finish L_{\lambda/\mu}\phi$, and $r_n>0$ 
otherwise. 
Now 
Corollary~\ref{nonzeroiff} yields $r_n=0$ and 
$I(d_n)=R$ for 
$1\le n\le \sum_{1+g\le t}k_t$ and for 
$n>\sum_{1\le t\le f}l_t$; 
hence (1) part (b) is true for $j\ge f+1$ and 
(1) part (c) is true for $j\le 0$. Thus for 
the rest of this proof we will also assume 
that $\sum_{1+g\le t}k_t < n \le \sum_{t\le f}l_t$. 
In particular, we also have $r_n > 0$, and by construction 
$I(d_n)\subseteq I_1(\phi)$.  

Now let $P$ be a prime ideal of $R$. 
If $I_1(\phi)\subseteq P$ then 
$P$ contains the radicals of both 
$I(d_n)$ and the corresponding 
ideal of minors of $\phi$ in each 
of the cases (a), (b), and (c). 
Therefore, we may assume without 
loss of generality that $I_1(\phi)$ is 
not contained in $P$, and it suffices 
to prove our 
statements after localization at $P$. 
In that case  $\phi_P=\phi_1\oplus\phi''$ 
with $\phi_1$ an isomorphism, 
and $\phi''\: F''\lra G''$ with 
$f''=\rank F''=f-1$ and $g''=\rank G''=g-1$. 
By Theorem~\ref{mapsplits} we 
have for each $n$ that 
$I(d_n)_P=I(d''_n)$ where $d''_n$ are 
the maps in $L_{\lambda/\mu}\phi''$. 
In view of Remark~\ref{R:threshold}(b) 
the desired conclusion 
is now immediate from our induction 
hypothesis applied to $L_{\lambda/\mu}\phi''$. 
\end{proof}

The following theorem generalizes  
\cite[Theorem 5.1.17]{ABW1982}, 
\cite[Theorem 2.1]{MS1996}, and 
\cite[Theorem 5.0.6]{A2017}.

\begin{thm} \label{acyclicity}
Let $\mu\subsetneq\lambda$, 
let $\phi: F\lra G$ be a 
map of finite free $R$-modules, 
let $f=\rank F$, \ 
let \ $g=\rank G$, \ and   
let \ $T=T_{\lambda/\mu}(f, g)$. 
The following are equivalent: 
\begin{enumerate}
\item 
The Schur complex 
$L_{\lambda/\mu}\phi$ is acyclic;

\item
$
\grade I_{f-j+1}(\phi) \ \geq \  
\sum_{i\le j}l_i 
$
\ for each \ $j\ge \max(1, \ f-T)$. 
\end{enumerate}
\end{thm}

\begin{proof}
When $T<0$ we have 
$L_{\lambda/\mu}\phi=0$ and both 
conditions of the theorem are trivially satisfied. When $T\ge 0$ 
the desired conclusion is 
immediate from Theorem~\ref{sqrt} and the Buchsbaum-Eisenbud acyclicity criterion.
\end{proof}

We say that $\lambda/\mu$ is a \emph{shift} 
by $(s,t)\in\mathbb N^2$ of a partition $\gamma$ if 
we have $\mu_j=t$ and $\lambda_j=t+\gamma_{j-s}$
for $s<j\le s+\tilde\gamma_1$,  
and $\mu_j=\lambda_j$ otherwise. 
The following acyclicity criterion holds 
in the generic case and 
generalizes in a different way 
\cite[Theorem 5.1.17]{ABW1982}.

\begin{theorem} \label{genericacyclicity}
Let $\mu\subsetneq\lambda$, 
let $\phi\: F\lra G$ be a generic map
with $f=\rank F\ge 1$ and $g=\rank G\ge 1$. 
The following are 
equivalent:
\begin{enumerate}
\item
The Schur complex $L_{\lambda/\mu}\phi$ is 
acyclic and nonzero. 

\item
%
\underline{Either}

\noindent
\quad{\rm (a)} \  
    $\lambda$ differs from        
    $\mu$ in at most $g-f+1$ columns; 

\noindent 
\underline{or}

\noindent
\quad{\rm (b)} \   
    $\lambda/\mu$ is the shift of a 
    partition $\gamma$ with \ 
    $g-f + 1< \gamma_1\le g$
    
\noindent
\quad\quad\quad    
    and \ 
    $\tilde\gamma_{\gamma_1}>\gamma_1-(g-f+1)$. 
\end{enumerate}
\end{theorem}

\begin{proof}
First we prove that (1) implies (2). 
The acyclicity and nonzero assumption, 
together with 
Proposition~\ref{P:nonzero-differential} 
and the fact that by construction 
$I_1(d_n)\subseteq I_1(\phi)$ for each $n$, 
imply that 
$L_{\lambda/\mu}(\phi)_0\ne 0$. 
Also, $T=T_{\lambda/\mu}(f,g)\ge 0$. 
Furthermore, 
by Theorem~\ref{eisenbudacyclicity} 
we must have 
$
\sqrt{I(d_1)}\subseteq \sqrt{I(d_2)}
\subseteq \dots \ .
$
Thus Theorem~\ref{sqrt}(a,c) yields 
$W=W(\lambda/\mu)\le g-T$, whence  $W= g-T$ by 
Remark~\ref{R:threshold}(c), and  
so $\nu''(\lambda,\mu,g-T)=\lambda$. 
Consider a box in the earliest possible row $i$ of 
$\nu'(\lambda,\mu,f-T-1)$ that is not in 
$\nu''(\lambda,\mu, W-1)$. Thus 
$\lambda_i-\mu_i=W$, and our box is in column 
$\lambda_i$. Also, our box is not in $\mu$. 
Now, we 
either have $f-T\le 1$, or there are at least  
$f-T-1\ge 1$ boxes below our box in column
$\lambda_i$. In the former case 
we obtain $f-T\le 1\le W = g-T$, hence $f\le g$, 
and our acyclicity 
assumption and Theorem~\ref{acyclicity} 
yield that 
$g-f+1=\grade I_f(\phi)\ge l_1$  
which implies (2) part (a). 
In the latter case for every box in 
column $\lambda_i$ and row $j>i$ we must have 
$\lambda_j=\lambda_i$ and so 
$\lambda_j-\mu_j\ge\lambda_i-\mu_i=W$, hence 
$\mu_j=\mu_i$. Therefore we have 
$l_t\ge W$ for $t=1,\dots,f-T$, which yields 
\[
l_1+\dots+ l_{f-T}\ge (f-T)W = (f-T)(g-T) = 
\grade I_{1+T}(\phi).  
\]
Our acyclicity assumption and 
Theorem~\ref{acyclicity} now 
imply $\grade I_{1+T}(\phi)=l_1+\dots +l_{f-T}$.  
But this equality is possible if and only if 
$l_t=W$ for $t=1,\dots,f-T$, 
hence in this case 
$\lambda/\mu$ is the shift by $(i,\mu_i)$ 
of a partition 
$\gamma$ with $g\ge \gamma_1=g-T> g -f + 1$ and 
$\tilde\gamma_{g-T}\ge f-T = \gamma_1 - g +f$, 
which yields 
(2) part (b). This completes the proof of 
$(1)\implies (2)$. 

Next, we show (2) implies (1). Suppose first 
that (2) part (a) holds. Then 
\[
g-f+1\ge l_1 \ge \dots \ge l_f, 
\]
and therefore 
$
\grade I_{f-k+1}(\phi)=k(g-f+k)\ge k(g-f+1)\ge 
l_1+\dots + l_k
$
for all $k\ge 1$. 
Thus Theorem~\ref{acyclicity}  
yields the desired acyclicity. Furthermore, 
$g-T\le W$ by Remark~\ref{R:threshold}. Thus 
$g-T\le W\le g -f +1\le g$ hence $T\ge 0$, yielding 
that $L_{\lambda/\mu}\phi$ is nonzero and 
completing the proof that (2) part (a) implies (1). 

Finally, suppose 
(2) part (b) holds, i.e. that 
$\lambda/\mu$ is the shift by $(s,t)$ of 
a partition $\gamma$ 
with $g\ge \gamma_1> g-f+1$ and
$\tilde\gamma_{\gamma_1}>\gamma_1 - (g-f+1)$. 
We want to show that 
$
\grade I_{f-k+1} (\phi) \geq l_1 + \cdots + l_k
$ 
for $k\geq f-T$ and use Theorem~\ref{acyclicity}. 

Note that in this case we must have $T=g-\gamma_1$. 
Indeed, 
\[
\nu''\bigl(\lambda,\mu, g - (g-\gamma_1)\bigr)=
\nu''(\lambda,\mu, \gamma_1)=\lambda
\]
hence it contains 
$
\nu'\bigl(\lambda,\mu, f - (g-\gamma_1)\bigr) = 
\nu'\bigl(\lambda,\mu, \gamma_1 - (g-f)\bigr). 
$ 
Also, row $s+1$ of 
$\nu''(\lambda,\mu,\gamma_1 - 1)$ has length 
$\mu_{s+1}+\gamma_1 -1$, while row $s+1$ of 
$\nu'\bigl(\lambda,\mu, \gamma_1 - (g-f+1)\bigr)$ 
has, in view of our assumption on 
$\tilde\gamma_{\gamma_1}$, 
length $\mu_{s+1}+\gamma_1=\lambda_{s+1}$. 
Thus $T=g-\gamma_1$ as claimed. In particular, 
$T\ge 0$ and so $L_{\lambda/\mu}\phi$ is nonzero. 

Now we have 
$f-T= f- (g-\gamma_1)=\gamma_1 - (g-f+1) + 1\ge 2$. 
Also,  
$l_i\leq \gamma_1$ for each $i$. Thus when 
$k\ge f-T$ we obtain $k\ge f-(g-\gamma_1)$ hence 
$1\le \gamma_1\le g-f+k$ and therefore 
for $f-T\le k\le f$ we get 
\[
l_1 + \cdots + l_{k} \leq k\gamma_1\leq 
k(g-f+k)=\grade I_{f-k+1}(\phi).
\] 
The desired acyclicity now follows from 
Theorem~\ref{acyclicity}.
\end{proof}

\begin{Example}
Let $\phi: F \rightarrow G$ be the generic map with $f=2$, $g=2$, $\lambda = (3,2,1)$, and $\mu = (2,2,1)$. Then $T=1$ and 
\begin{align*}
&\grade I_2(\phi) = 1 \geq l_1 = 1 \\
&\grade I_1(\phi) = 4 \geq l_1 + l_2 = 1 \\
&\grade I_{f-k+1}(\phi) = \infty \geq 1 \text{ for all } k \geq 3 
\end{align*}
so that $L_{\lambda/\mu}\phi$ is acyclic by Theorem~\ref{genericacyclicity}.
\end{Example}

\begin{Example}
Let $\phi:F \rightarrow G$ be the generic map with $f=4$, $g=2$, and $\lambda=(1,1,1)$. Then $\gamma_1=1$ and $\tilde{\gamma}_{\gamma_1} = 3$ so that condition (2)(b) is satisfied in Theorem \ref{genericacyclicity} and it follows that $L_{\lambda/\mu}\phi$ is acyclic. 
\end{Example}

\begin{Example}
Let $\phi: F \rightarrow G$ be generic with $f=2, g=4, \lambda = (4,4,1), \text{ and } \mu = (0)$. Since $3< \gamma_1 \leq 4$ and $\tilde{\gamma}_{\gamma_1}>1$, Theorem \ref{genericacyclicity} insures that the complex $L_{\lambda/\mu}\phi$ is acyclic.
\end{Example}

\section{Torsion freeness of Schur modules}
\label{S:torsion-freeness}

Everything is now in place to begin our study 
of the torsion freeness of Schur modules. 
First, we introduce a chain complex with good 
rigidity properties.

\begin{definition}
Let 
$
F \xrightarrow{\phi} G 
$
be a map of free $R$-modules with 
$f=\rank F<\rank G$. 
We write $\tilde{L}_{\lambda/\mu}\phi$ 
for the chain complex
$$
0 \rightarrow C_n \rightarrow C_{n-1} 
\rightarrow \cdots \rightarrow L_{\lambda/\mu}G
\xrightarrow{L_{\lambda/\mu}(\epsilon^*)} 
L_{\lambda/\mu}(\wedge^{f+1}G\otimes\wedge^f F^*)
\rightarrow 0,
$$
where $\epsilon=\epsilon(\phi)$ is the map 
from \eqref{E:epsilon-def}, and the piece 
$$
0 \rightarrow C_n \rightarrow C_{n-1} \rightarrow \cdots \rightarrow L_{\lambda/\mu}G 
$$
is just the Schur complex $L_{\lambda/\mu}\phi$ 
shifted so that 
$L_{\lambda/\mu}G$ is in homological degree $1$.  
\end{definition}

\begin{remarks} 
(a) 
Since by Remark~\ref{R:epsilon-functorial}(a) 
the map $\epsilon^*$ factors through the 
canonical projection 
$\pi\: G \longrightarrow M= G/\im\phi$, it 
follows from Proposition~\ref{P:Schur-homology} 
that $\tilde L_{\lambda/\mu}\phi$ is indeed 
a chain complex. 

(b) 
It is straightforward from the definition that 
the chain complex 
$\tilde{L}_{\lambda/\mu}\phi$ commutes with 
base change, and therefore, in view of 
Remark~\ref{R:epsilon-functorial}(b), it 
induces a layered functor 
\smash[t]{
$
\tilde{L}_{\lambda/\mu}\: 
\widetilde{\mathcal{M}}_{g,f} \rightarrow \text{\bf Comp}.
$
}
\end{remarks}

The first homology of $\tilde L_{\lambda/\mu}\phi$ 
has the following key property. 

\begin{lemma} \label{exactatllambdag}
Let $\phi: F \rightarrow G$ be an 
injective map of free 
$R$-modules with 
$f = \operatorname{rank}F < \rank G$, 
let $M=\coker\phi$, and let \ 
$
\overline{\epsilon^*}\: 
M \longrightarrow 
\wedge^{f+1}G\otimes\wedge^f F^* 
$
be the induced by \ $\epsilon^*$ map.  
The following are equivalent: 
\begin{enumerate}
\item 
$L_{\lambda/\mu} M$ is torsion free;  

\item 
The map \ 
$
L_{\lambda/\mu}
\bigl(\overline{\epsilon^*}\bigr)  
\: \ 
L_{\lambda/\mu}M \longrightarrow 
L_{\lambda/\mu} 
(\wedge^{f+1}G\otimes \wedge^f F^*)
$ 
\ is injective.

\item
$\HH_1(\tilde{L}_{\lambda/\mu}\phi)=0$.
\end{enumerate} 
\end{lemma}

\begin{proof}
Let $\pi\: G \longrightarrow M=G/\im\phi$ be 
the canonical projection, so that 
$\epsilon^*=\overline{\epsilon^*}\pi$. The 
equivalence of (2) and (3) is immediate from 
Proposition~\ref{P:Schur-homology}. 
As (2) implies (1) trivially, we proceed to  
show that (1) implies (2). Thus, we may 
assume all non-zerodivisors are units, hence 
$M$ is free and $\phi$ splits. 
We may therefore further assume,
as in the proof of Lemma~\ref{exactatg},
that $R$ is a field. But then by 
Lemma~\ref{exactatg} we have that 
$\overline{\epsilon^*}$ 
is injective, hence split. 
Since split injections  
are functorial, 
$
L_{\lambda/\mu}
\bigl(\overline{\epsilon^*}\bigr)
$
must also be split injective. 
\end{proof}

We also need the following immediate 
consequence of Theorem~\ref{acyclicity}. 

\begin{corollary} \label{C:torsion-freeness}
Let $\mu\subsetneq\lambda$, let $\phi: F\lra G$ be 
a map of finite free $R$-modules, 
let $f=\rank F$, \ 
let \ $g=\rank G$, \ 
let \ $T=T_{\lambda/\mu}(f, g)$, \ and 
let $M=\coker \phi$. 
The following are equivalent: 
\begin{enumerate}
\item 
The complex 
$L_{\lambda/\mu}(\phi)$ is acyclic and 
$L_{\lambda/\mu}M$ is torsion-free. 

\item
$
\grade I_{f-j+1}(\phi) \ \geq \  
1 + \sum_{i\le j}l_i 
$
\ for each \ $j\ge \max(1, \ f-T)$. 
\end{enumerate}
\end{corollary}

This allows us to obtain the following result 
in the generic case. 

\begin{theorem}\label{T:generic-torsionfreeness}
Let $\mu\subsetneq\lambda$, and let 
$\phi\:F \lra G$ be 
generic with $f=\rank F\ge 1$ and 
$g=\rank G\ge 1$.
Let $M=\coker\phi$. 
The following are equivalent: 
\begin{enumerate}
    \item 
    The complex $L_{\lambda/\mu}\phi$
    is acyclic, and $L_{\lambda/\mu}M$ is 
    nonzero and torsion-free.  
    
    \item
    The partition $\lambda$ differs from $\mu$ 
    in at most $g-f$ columns. 
\end{enumerate}
\end{theorem}

\begin{proof}
We show that (2) implies (1). Indeed, by (2) 
we have $g-f\ge l_1 \ge \dots \ge l_f$, 
therefore for all $1\le k\le\min(f,g)$ we have 
\[
\grade I_{f-k+1}(\phi) = k(g-f+k) =  
k^2 + k(g-f) \ge 1 + k(g-f)\ge 1 + \sum_{i\le k}l_i,  
\]
and the desired acyclicity and torsion 
freeness follow from 
Corollary~\ref{C:torsion-freeness}. Since 
by construction of the differentials $d_i$ 
of $L_{\lambda/\mu}\phi$ one has 
$
\im d_i\subseteq
I_1(\phi)(L_{\lambda/\mu}\phi)_{i-1}
$
for each $i$, acyclicity together with the 
fact that $\phi$ is generic imply that 
$L_{\lambda/\mu}G$, and hence $L_{\lambda/\mu}M$, 
are nonzero.

For (1) implies (2), we must rule out (b) in Theorem \ref{genericacyclicity} (2). Assume that $\gamma_1 = g-f+s$ where $s \geq 1$ and that $\tilde{\gamma}_{\gamma_1} = \gamma_1 - (g-f+1) = (g-f+s) - (g-f+1) = s-1$. We note that $\nu''$ starts getting smaller when $g-t < \gamma_1 = g-f+s \implies f-s<t$. In particular $\nu_1$ changes from $\lambda_1$ to $\lambda_1 - 1$ when $t = f-s+1$. Also, $\nu_1' = \lambda_1$ when $f-t < \tilde{\gamma}_{\gamma_1} \implies f- \tilde{\gamma}_{\gamma_1} < t$. 

By assumption, we have $\tilde{\gamma}_{\gamma_1} \geq s$ so that $f- \tilde{\gamma}_{\gamma_1} \leq f-s$. It follows that $f-T=f-(f-s)=s$. But then $\grade I_{f-s+1}(\phi) = (f-f+s-1+1)(g-f+s-1+1) = (s)(g-f+s) \not\geq s\gamma_1+1 = s(g-f+s)+1$. Thus we have ruled out (b) and the case in (a) when $\lambda$ differs from $\mu$ in $g-f+1$ columns and the result follows. 
\end{proof}

We will need the following key consequence of 
Theorem~\ref{T:generic-torsionfreeness}. 

\begin{corollary}\label{C:tilde-rigidity} 
Suppose $\mu\subsetneq\lambda$. 
The functor 
$
\tilde{L}_{\lambda/\mu}\: 
\widetilde{\mathcal M}_{g,f}\longrightarrow
\text{\bf Comp}
$ 
is rigid when $\lambda$ differs from $\mu$ in at most $g-f$ columns.
\end{corollary}

\begin{proof}
Since $\mu\ne \lambda$ and they differ in 
at most $g-f$ columns, we must have $f<g$. 
By the rigidity criterion, 
we only need to show that 
\smash[t]{
$\tilde L_{\lambda/\mu}\phi$
}
is acyclic in the 
generic case. Note that if $f=0$ then 
$\tilde L_{\lambda/\mu}\phi$ has the form 
\smash[t]{
$
0\lra L_{\lambda/\mu}G \xrightarrow{\ 1 \ } 
      L_{\lambda/\mu}G \lra 0
$
}
hence is acyclic. Thus we may assume $f\ge 1$. 
But then by 
Theorem~\ref{T:generic-torsionfreeness} 
the complex $L_{\lambda/\mu}\phi$
is acyclic, and 
$L_{\lambda/\mu}M$ is torsion free.
The desired acuclicity is now 
immediate from 
Lemma~\ref{exactatllambdag}. 
\end{proof}

The following is the ``local'' version 
of our second main result:

\begin{thm}\label{mainschur}
Let $\phi: F \rightarrow G$ be an 
injective map of free $R$-modules 
with $f = \rank F< \rank G = g$. Let 
$M$ be the cokernel of $\phi$. 
Let $\mu\subsetneq\lambda$ be partitions  
such that $\lambda$ differs from $\mu$ in at most 
$g-f$ columns, and let $T=T_{\lambda/\mu}(f,g)$. 

Then $L_{\lambda/\mu}M$ is torsion free if 
and only if \ 
$
\grade I_{f-k+1}(\phi) \geq 1 + 
\sum_{i\le k}l_i
$ 
\ for all \ $k\ge \max(1, \ f-T)$.
\end{thm}

\begin{proof}
If 
$
\operatorname{grade} I_{f-k+1}(\phi) \geq
1+\sum_{i\le k}l_i
$ 
for all 
$k\ge\max(1, f-T)$, then 
by Corollary~\ref{C:torsion-freeness}  
the complex 
$L_{\lambda/\mu}\phi$ is acyclic and 
$L_{\lambda/\mu}M$ is torsion free. 

If $L_{\lambda/\mu}M$ is torsion free, then  
by Lemma~\ref{exactatllambdag} we have  
$\HH_1(\tilde{L}_{\lambda/\mu}\phi) = 0$. 
Since $\tilde{L}_{\lambda/\mu}$ is rigid 
by Corollary~\ref{C:tilde-rigidity}, the 
complex $L_{\lambda/\mu}\phi$ is acyclic and the
result follows from 
Corollary~\ref{C:torsion-freeness}. 
\end{proof}

Now we are ready to prove our second main 
result. 

\begin{theorem}\label{T:second-main} 
Let $M$ be a finitely generated $R$-module 
with $\operatorname{pd}_R M\leq 1$ and 
$\rank M = r$. Let $\mu\subsetneq \lambda$ be 
partitions such that $\lambda$ differs from $\mu$ 
in at most $r$ columns, and let 
$T = T_{\lambda/\mu}(0,r)$. 

Then $L_{\lambda/\mu}M$ is torsion free 
if and only if \ 
$
\grade \operatorname{Fitt}_{r+k-1}(M) \geq 
1 + \sum_{i\le k}l_i
$ 
\ for all \ $k\ge \max(1, -T)$.
\end{theorem}

\begin{proof}
Since $L_{\lambda/\mu}M$ is finitely generated 
and has rank, 
the module $L_{\lambda/\mu}M$ is torsion free if 
and only if $L_{\lambda/\mu}M_{\mathfrak{p}}$ is 
torsion free for all $\mathfrak{p} \in \Spec{R}$. 
Thus we may assume that 
$R$ is local, in particular $M$ has a 
minimal free resolution $0\lra F\lra G\lra 0$.  
Let $f=\rank F$ and $g=\rank G$. 
Since $\mu\ne\lambda$ and $\lambda$ differs from 
$\mu$ in at most $r$ columns we get $1\le r=g-f$,
hence $f<g$. Also, by 
Remark~\ref{R:threshold}(b) we have 
$T_{\lambda/\mu}(f,g) - f = T_{\lambda/\mu}(0,r)$, 
and the desired result follows from 
Theorem~\ref{mainschur}.
\end{proof}




\end{document}